\documentclass[12]{article}

\usepackage{authblk}
\usepackage{amsmath} % provides many mathematical environments & tools
\usepackage{amsfonts}
\usepackage{graphicx}
\usepackage{float}

\usepackage{setspace}

\newcommand{\norm}[1]{\left\lVert#1\right\rVert}
\newcommand{\myat}[1]{\left. #1 \right|}
\newcommand{\mydelta}[1]{\, \delta #1 \,}

\begin{document}

\title{Iterative PDE-constrained optimization for seismic full-waveform inversion}

\author[1]{Mikhail Malovichko}
\author[1]{Akbar Orazbayev}
\author[1]{Nikolay Khokhlov}

\affil[1]{Moscow Institute of Physics and Technology}

\date{\today}
\maketitle

\begin{abstract}

This paper presents a novel numerical method for the Newton seismic full-waveform inversion (FWI).
The method is based on the full-space approach, where the state, adjoint state, and control variables are optimized simultaneously.
Each Newton step is formulated as a PDE-constrained optimization problem, which is cast in the form of 
the Karush-Kuhn-Tucker (KKT) system of linear algebraic equitations.
The KKT system is solved inexactly with a preconditioned Krylov solver.
We introduced two preconditioners: the one based on the block-triangular factorization and its variant with an inexact block solver.
The method was benchmarked against the standard truncated Newton FWI scheme on a part of the Marmousi velocity model.
The algorithm demonstrated a considerable runtime reduction compared to the standard FWI. 
Moreover, the presented approach has a great potential for further acceleration.
The central result of this paper is that it establishes the feasibility of Newton-type optimization of the KKT system in application to the seismic FWI.

\end{abstract}

\section{Introduction}

The seismic full-waveform inversion (FWI) is a coefficient inverse problem aiming to estimate the subsurface distribution of material properties based on recorded seismic data.
The FWI has emerged almost four decades ago \cite{Lailly1983,Tarantola1984,Mora1987}, 
and has evolved into a powerful method for reconstruction of the subsurface properties, 
see reviews \cite{Virieux2009,Tromp2020}. 
 
The FWI is formulated as a Tikhonov-like minimization procedure,
\begin{equation}
	\label{eq:basic_tikhonov}
	\underset{\gamma}{\text{minimize}} \quad \norm{d-\mathcal{F}(\gamma)}^2 + \mathcal{R}(\gamma),
\end{equation}
where $\gamma$ is a function of material parameters (say, the squared P-wave slowness),
$d$ is a vector of observed data, $\mathcal{F}$ is a nonlinear forward-problem operator, and $\mathcal{R}$ is the stabilizing term.
We are interested in the frequency-domain acoustic FWI and, thus, operator $\mathcal{F}$ implies the solution of Helmholtz's equation.
Many variants of basic formulation \eqref{eq:basic_tikhonov} exist: $\mathcal{R}$ can take various forms or even be omitted, the residual can be weighted, and so on.

Optimization \eqref{eq:basic_tikhonov} is usually performed by either quasi-Newton or Newton methods.
The quasi-Newton methods dominate industrial applications due to their relative numerical efficiency
\cite{Plessix2009,Warner2013,Operto2015,Operto2018,Fichtner2009} (among many others).
These methods account for nonlinearity partially, usually by approximating the Hessian with a diagonal matrix. 
Typical representatives of this class are the nonlinear conjugate gradient (NLCG) method, where a diagonal matrix scales the gradients \cite{Shin2001,Mulder2008}, and L-BFGS \cite{Brossier2009}. 
Though computationally tractable, the quasi-Newtonian methods may suffer from convergence stagnation. 

The Newton methods achieve faster and more reliable convergence,
which catch up with the nonlinearity through quadratic approximation of the target function, see 
\cite{Pratt1998,Epanomeritakis2008,Metivier2017}.
The often-used Gauss-Newton variant computes the Hessian matrix approximately by linearizing the forward-problem operator.
The main disadvantage of the Newton methods applied to \eqref{eq:basic_tikhonov} is an enormous computational burden stemming from the need to solve a large dense system of linear equations at each Newton iteration. 
The system is solved iteratively, usually with the conjugate gradients (CG). 
The matrix is ill-conditioned, and therefore each Newton (external, nonlinear) iteration may require hundreds of internal (linear) iterations. 
In addition, no efficient preconditioner is known to improve this matrix's conditioning.
For these reasons, numerical implementations of the Newton FWI are far less common than the quasi-Newton ones.

It is natural to ask if we can obtain the fast Newton-like convergence with a less resource-demanding procedure.
The ultimate goal of this paper is to demonstrate that the answer is yes.

Optimization \eqref{eq:basic_tikhonov} is known as \textit{the reduced-space approach} in the optimization community.
Its alternative is the \textit{the full-space approach}, also known as \textit{all-at-once}.
In the reduced-space approach, only material parameters $\gamma$ are optimized, whereas the connection between the wavefield in the medium and $\gamma$ is prescribed by the forward-problem operator $\mathcal{F}$.
In the full-space approach, one simultaneously optimizes the control, state, and adjoint-state variables.
Upon formulating the necessary conditions and discretizing the PDE, one obtains a nonlinear system of algebraic equations, known as the Karush-Kuhn-Tucker (KKT) system, which is sparse. 
Its sparsity is a clear advantage over formulation \eqref{eq:basic_tikhonov}.
The main challenge here is that the problem becomes many times bigger than the reduced-space approach because the state and adjoint-state variables for all sources are included in optimization.

The full-space approach has been introduced to geophysics in the seminal paper \cite{Haber2000}, followed by \cite{Haber2001,Haber2004}.
Those authors cast the geophysical electromagnetic inversion as a full-space optimization problem, 
derived the KKT system and applied the Newton algorithm to it.
Despite spectacular results, the full-space approach had been considered impractical for years due to the enormous size of the linear system, 
arising on each Newton step.
Regarding FWI, the KKT optimality conditions were derived 
in several papers \cite{Epanomeritakis2008,Metivier2017}, 
but were viewed as a theoretical tool rather than a basis for numerical implementation.

Significant progress has been made in this direction during the last decade.
Papers \cite{Abubakar2009,vanLeeuwen2013,vanLeeuwen2016} report the application of the \textit{penalty method} to the FWI, marking an important step toward the full-space approach.
The penalty method has some theoretical drawbacks connected to the choice of the penalty parameter \cite{Nocedal} (but see recent developments \cite{Rizzuti2021}).
The \textit{augmented Lagrangian} (AL) method is a viable alternative to the penalty method.
A variant of it, known as the \textit{alternating-direction method of multiplies} (ADMM) \cite{Boyd2010}, was applied to the FWI with great success, see \cite{Aghamiry2019, Aghazade2022}. 
In essence,  the ADMM tackles the full-space optimization by making alternating steps in the state, adjoint-state, and control state spaces.
Still, it is generally accepted that the main weakness of the ADMM compared to the Newton iterations is its slow convergence
\cite{Boyd2010}.

This paper studies the full-space optimization
based on applying a Newton solver to the original KKT matrix. 
The crucial component of our approach is a suitable preconditioner based on a block-triangular factorization of the permuted KKT matrix
combined with approximations of individual blocks.
Similar ideas were exploited optimization with Maxwell's equations \cite{Haber2001} and the Stokes equation \cite{Biros2005}; thus, this paper extends this framework to Helmholtz's equation.
We benchmarked our method against the standard reduced-space Gauss-Newton FWI with the internal CG solver.
In the numerical experiments conducted the Marmousi model, 
we observed the 6$\times$ speedup.
The method has the potential for further acceleration. 
In particular, it will benefit from any method designed to approximate the reduced-space Hessian, such as L-BFGS.
The central result of this paper is that it establishes the feasibility of Newton-type optimization of the KKT system in application to seismic FWI.

The paper is organized as follows.
In Section~2, we formulate the forward problem.
The reduced-space Gauss-Newton inversion is considered in Section~3.
Sections~4 and 5 are dedicated to the full-space approach and its iterative solution.
Numerical experiments are presented in Section~6.
Concluding remarks are given in Section~7.

\section{The forward problem}

\begin{figure}[h!]
	\centering
	\includegraphics[width=8cm]{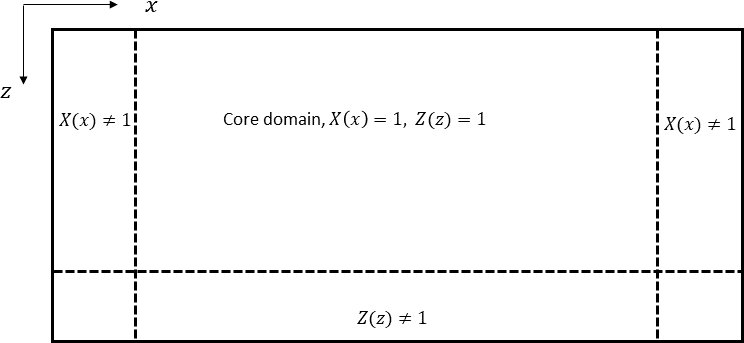}
	\caption{Schematic representation of the computational domain with a PML region.
	In the PML region $X(x) \neq 1$,  $Z(z) \neq 1$.
	}
	\label{fig:PML}
\end{figure}

In this section, we briefly review the forward problem.
Let $\Omega \subset \mathbb{R}^2$ be a bounded domain
with boundary $\partial \Omega$.
Consider the 2D Helmholtz's equation in the following form,
\begin{equation}
	\label{eq:helmholtz_preliminary}
	\begin{aligned}
	-Z \frac{\partial}{\partial x} \left( \frac{1}{X} \frac{\partial U}{\partial x} \right)
	-X \frac{\partial}{\partial z} \left( \frac{1}{Z} \frac{\partial U}{\partial z} \right)
	- 
	\omega^2 \gamma XZ U 
	= 
	XZF, 
	\text{ in } \Omega, \\
	U = 0, \text{ on } \partial \Omega.
	\end{aligned}
\end{equation}
Here $F(x,z)$ are the right-hand side, 
$U(x,z,\omega)$ are the acoustic field (monochromatic pressure), 
$\omega$ is the circular frequency,
$\gamma(x,z)$ is the squared slowness, $\gamma=c^{-2}$, where $c(x,z)$ is the speed of sound.
Complex-valued functions of one variable $X(x)$ and $Z(z)$ implement the perfectly-matched layers (PML).
They deviate from 1.0 only in the damping region on the sides of the computational domain to suppress the side reflections, 
see Figure~\ref{fig:PML}.
On the PML implementation, see \cite{Yavich2021} and references therein.

Let us rewrite \eqref{eq:helmholtz_preliminary} with a more compact notation as follows,
\begin{equation}
	\label{eq:helmholtz_preliminary2}
	\begin{aligned}
	\mathcal{D} U = F
	\text{ in } \Omega, \\
	U = 0, \text{ on } \partial \Omega.
	\end{aligned}
\end{equation}
where 
\begin{equation}
	\mathcal{D}(\cdot) := 
	-Z \frac{\partial}{\partial x} \frac{1}{X} \frac{\partial}{\partial x}(\cdot)
	-X \frac{\partial}{\partial z} \frac{1}{Z} \frac{\partial}{\partial z}(\cdot)
	- 
	\omega^2 \gamma XZ (\cdot) \,.
\end{equation}
Note that we replaced $XZF$ with $F$ in \eqref{eq:helmholtz_preliminary2} because we assume that 
$F$ is supported only in the core domain.

In a seismic survey, the acoustic field is excited by many monochromatic sources independent of each other.
Thus, we have a series of problems,
\begin{equation}
	\label{eq:helmholtz_j}
	\begin{aligned}
		\mathcal{D} U_k = F_k, \text{ in } \Omega, \\
		U_k = 0, \text{ on } \partial \Omega,
	\end{aligned}
\end{equation}
where index $k$ designates the source index.
Let us combine solutions $U_k$ and right-hand sides $F_k$ into vector variables to make the formulas even more compact.
Starting from now on, $U$ is a direct product of all $U_k$, $F$ is a direct product of $F_k$:
\begin{equation}
	U := \prod_{k=0}^{K-1} U_k, \quad F := \prod_{k=0}^{K-1} F_k, \quad k=0..K-1.
\end{equation}
With these notations, a set of $K$ problems \eqref{eq:helmholtz_j} can be compactly written as
\begin{equation}
	\label{eq:helmholtz}
	\begin{aligned}
		\mathcal{D} U = F, \text{ in } \Omega, \\
		U=0 , \text{ on } \partial \Omega.
	\end{aligned}
\end{equation}
Now we introduce the observation operator $\mathcal{Q}$ which takes the acoustic field and returns a vector of data $d$.
We define $\mathcal{Q}$ through the convolution with the Dirac $\delta$-function (although any linear functional can be used instead) as follows,
\begin{equation}
		d = \mathcal{Q} U := \left\{ \int_{\Omega} \delta( x_{kj}-x) U_k \, dx \right\},
		\quad k=0,1,..K-1, \quad j=0,1,..N_k-1
\end{equation}
where $x_{kj}$ is the position of the $j$-th receiver of the $k$-th source, 
$N_k$ is the number of receivers for the $k$-th source,
$d \in \mathbb{C}^N$, $N=\sum_{k=0}^{K-1}N_k$ is the total number of receivers for all sources.

In what follows, we will use the discrete form of the forward problem.
The finite-difference (FD) discretization is assumed
because it is used in our numerical experiments.
However, we might as well have performed discretization by the finite-element of spectral-element methods.

Let $u_k$, $f_k$, and $s$ be finite-difference approximations of $U_k$, $F_k$, and $\gamma$, respectively.
We introduce the following notation,
\begin{equation}
	u = \prod_{k=0}^{K-1} u_k, \quad f = \prod_{k=0}^{K-1} f_k, \quad k=0..K-1.
\end{equation}
A set of discrete forward problems reads
\begin{equation}
	\begin{split}
		\label{eq:Au=f}
		A u = f,
	\end{split}
\end{equation}
where $A$ is a block-diagonal matrix.
Each diagonal block $A_j$ of matrix $A$ is a sparse complex non-Hermitian matrix.
If the PMLs are not used, then $A_j$ turns into the real symmetric indefinite matrix,
\begin{equation}
		-\Delta- \omega^2 M,
\end{equation}
where $\Delta$ is a discrete Laplacian, $M$ is a diagonal mass matrix, containing values of $s$ on its diagonal $M= I s$ (it can be non-diagonal in some FD schemes \cite{Hustedt2004}).

%In our numerical experiments, which are two-dimensional, we solve each equation $A_j u_j=f_j$ by the LU factorization.

\section{The reduced-space Gauss-Newton inversion}

In this section, we review the reduced-space Gauss-Newton inversion (RSGN).
Since we will eventually end up with the FD approximation, and thus discrete variables will be grid functions, it is easier to formulate the inverse problem in the discrete form (the \textit{discretize-then-optimize} approach).
In the finite-element framework, a more convenient way would be to formulate the optimality conditions in function spaces and then move on to discretization (the \textit{optimize-then-discretize} approach), see \cite{Malovichko}.

Let $S$ be a nonlinear operator that maps a given $s$ to $u$,
\begin{equation}
	u = S(s).
\end{equation}
Let $Q$ be a finite-dimensional operator that maps a given discrete acoustic field $u$ to data $d$,
\begin{equation}
	d = Q u.
\end{equation}
Discrete coefficient $s$ is iterated as $s_{n+1} = s_n + \mydelta{s}$, where $\mydelta{s}$, known as the model update.
On each iteration, $s$ is a solution to the following minimization problem,
\begin{equation}
		\label{eq:tikhonov_discrete}
		\underset{\delta s}{\text{minimize}} \quad
		\frac{1}{2} 
		\norm{ W(r - Q G \mydelta{s})}^2_2 + 
		\frac{\varepsilon}{2} \norm{ s_n+\mydelta{s} }_2^2.
\end{equation}
Here $W$ is a real-valued diagonal matrix of weights designed primarily to balance signal attenuation with offset from the source.
Vector $r$ is the data residual,
\begin{equation}
	r = d-Q S(s_n),
\end{equation}
$G$ is a matrix of the Frech{\'e}t operator of $S$ at point $s=s_n$, e.i.
\begin{equation}
	\label{eq:G}
	G = \myat{ \frac{\partial S}{\partial s} }_{s=s_n}.
\end{equation}
%
%Note that in the Gauss-Newton method, the equality in \eqref{eq:G} is approximate because the second-order derivatives are dropped.
%
Solution to \eqref{eq:tikhonov_discrete} satisfies 
\begin{equation}
	\label{eq:normal_system_discrete}
	\left(
		G^* Q^* W^T W Q G + L
	\right) \mydelta{s}
	=
	W^T G^* Q^* r - L s_n,
\end{equation}
where $L=\varepsilon I$.
The following formula defines matrix $G$,
\begin{equation}
	\label{eq:S}
	G = A^{-1} P,
\end{equation}
where $P$ is a diagonal matrix containing forward solution $u_n$ at $s=s_n$,
\begin{equation}
	\label{eq:P}
	P := \omega^2 I u_n.
\end{equation}
To verify \eqref{eq:S},\eqref{eq:P}, one can subtract the two discrete Helmholtz's equations
\begin{equation}
	\begin{aligned}
		\left(-\Delta-\omega^2 I \, (s_n+\mydelta{s}) \right) (u_n+\delta u) = f, \\
		\left(-\Delta-\omega^2 I \, s_n \right) u_n = f, \\
	\end{aligned}
\end{equation}
although the derivation can be conducted directly in the Hilbert spaces.
%see also \cite{Plessix2006}
In geophysics, the matrix 
\begin{equation}
	\label{eq:jacobian}
	J := Q A^{-1} P,
\end{equation}
is known as the Jacobian,
the matrix 
\begin{equation}
	\label{eq:Hessian}
	H := J^*W^T W J + L
\end{equation}
is the (regularized) Hessian, 
and vector
\begin{equation}
	g := W^T J^* r - L s_n
\end{equation}
is called the gradient.
Using this definitions, we rewrite \eqref{eq:normal_system_discrete}
as 
\begin{equation}
	\label{eq:normal_system_discrete2}
	H \mydelta{s}
	=
	g.
\end{equation}

For any $\varepsilon>0$ matrix $H$ is positive-definite, and thus 
\eqref{eq:normal_system_discrete2} has a unique solution.
It should be noted that, even if $\varepsilon =0$ (singular $H$), then a Krylov solver will deliver a meaningful estimate of $\mydelta{s}$, that is, a projection to a corresponding Krylov subspace.

The system matrix in \eqref{eq:normal_system_discrete2} is dense and large. 
In real-life 3D problems, the system \eqref{eq:normal_system_discrete2} can be solved only iteratively, with the CG being a natural choice.
We will abbreviate this algorithm as RSGN-CG.
The CG multiplies the system matrix $J^* J$ by a vector once per iteration. 
Because of \eqref{eq:jacobian} it implies two linear solves per CG iteration - the forward simulation with $A$ and the adjoint simulation with $A^*$.
We can rewrite \eqref{eq:normal_system_discrete2} as 
\begin{equation}
	\label{eq:normal_system_discrete_sum}
	\left( \sum_{j=0}^{K-1} J^*_j W_j^T W_j J_j + L  \right) \mydelta{s}
	= \sum_{j=0}^{K-1} W_j^T J^*_j r_j - L s_n.
\end{equation}
The forward and adjoint simulations for different sources can be performed in parallel 
but these results must be combined at each CG iteration.

However, the most critical problem is a very slow convergence of the CG. 
The system matrix \eqref{eq:normal_system_discrete2} has a very high condition number for all meaningful values of $\varepsilon$.
Thus, tens or hundreds of forward simulations (for all sources) may be needed to generate a good update $\mydelta{s}$.
In addition, preconditioning of a regularized normal system of linear equitations is notoriously difficult.

%%%%%%%%%%%%%%%%%%%%%%%%%%%%%%%%%%%%%%%%%%%%%%%%%%%%%
%%%%%%%%%%%%%%%%%%%%%%%%%%%%%%%%%%%%%%%%%%%%%%%%%%%%%
\section{The full-space Gauss-Newton inversion}

In this section, we describe the full-space Gauss-Newton (FSGN) inversion.
We start by noting that 
the discrete optimization problem \eqref{eq:tikhonov_discrete} is equivalent to the following 
constrained optimization problem,
\begin{equation}
	\label{eq:constrained_optimization}
	\begin{split}
		\underset{\mydelta{s}}{\text{minimize}} \quad
		\Psi :=
		\frac{1}{2} 
		\norm{W \left( r - Q \mydelta{u} \right)}^2_2 + 
		\frac{\varepsilon}{2} \norm{ s_n+\mydelta{s} }_2^2,
		\\
		\text{ s.t. } \quad		A \mydelta{u} = P \, \mydelta{s}.
	\end{split}
\end{equation}
Here $n$ is the Newton iteration number, $\mydelta{u}$ is the update to $u_n$, e.i. $u_{n+1} \approx u_n + \mydelta{u}$.
The target functional $\Psi$ is quadratic and bounded away from zero for any $\varepsilon >0$.
The constraints are linear.
For such problems, known as linear-quadratic, the existence 
of unique state and control variables are readily established, for example, \cite{Hinze2009}[Thm 1.43].
Moreover, since we start from the discrete formulation, the solvability is proved simply by 
establishing the equivalence between the matrix form and the normal system 
\eqref{eq:normal_system_discrete2}, see below.

Let us form the Lagrangian,
\begin{equation}
	\mathcal{L} := \Psi + \lambda^* (A \mydelta{u} - P \, \mydelta{s}),
\end{equation}
where $\lambda$ is the adjoint-state variable (a complex-valued grid function).
By applying the optimality conditions \cite{Hinze2009}, we obtain the KKT system,
\begin{equation}
	\begin{split}
	-Q^* W^T W(r-Q \mydelta{u}) + A^* \lambda = 0, \\
	L (s_n + \mydelta{s}) - P^* \lambda = 0, \\
	A \mydelta{u} - P \mydelta{s} = 0.
	\end{split}
\end{equation}
This system can be rewritten in the matrix form,
\begin{equation}
\label{eq:KKT}
	\underbrace{
		\left[
		\begin{array}{ccc}
			 F & O   &  A^* \\
			 O & L   & -P^* \\
			 A &-P   & O    \\
		\end{array}
		\right]
	}_{\mathcal{M}}		
	\underbrace{
		\left[
		\begin{array}{c}
			\mydelta{u}  \\
			\mydelta{s}  \\ 
			\lambda
		\end{array}
		\right]
	}_{\xi}		
	= 
	\underbrace{
		\left[
		\begin{array}{c}
			Q^* W^T W r \\
			-L s_n \\ 
			0
		\end{array}
		\right],
	}_{b}		
\end{equation}
where $F = Q^* W^T W Q$.

Lets us prove that \eqref{eq:KKT} is equivalent to \eqref{eq:normal_system_discrete2} for it establishes solvability of the KKT system.
Using the third equation of \eqref{eq:KKT}, we compute $\mydelta{u}$,
\begin{equation}
	\label{eq:tmp1}
	\mydelta{u}  = A^{-1} P \mydelta{s}.
\end{equation}
Using this result, we compute $\lambda$ from the first equation,
\begin{equation}
	\label{eq:tmp2}
	\lambda  = A^{-*} \left(Q^* r - Q^* Q A^{-1} P \mydelta{s} \right),
\end{equation}
where notation $A^{-*}$ means $(A^*)^{-1}$.
Now we substitute \eqref{eq:tmp1} and \eqref{eq:tmp2} into the second equation of \eqref{eq:KKT} and get
\begin{equation}
	\label{eq:normal_system_discrete3}
	\left(  P^* A^{-*} Q^* Q A^{-1} P + L  \right) \mydelta{s} = P^* A^{-*} Q^* r  -L s_n.
\end{equation}
From \eqref{eq:jacobian} we see that \eqref{eq:normal_system_discrete3} is equivalent to
\eqref{eq:normal_system_discrete2}.

From the preceding discussion, it is clear that the model update $\mydelta{s}$ computed from either the KKT system \eqref{eq:KKT} or the normal system \eqref{eq:normal_system_discrete2} is identical, provided that the linear systems are solved exactly.
The difference between the two approaches lies in the structure of the system matrix.
The system matrix in \eqref{eq:normal_system_discrete2} is complex-valued, dense, Hermitian positive-definite, and relatively small.
For small- to medium-size 3D problems with rigorously selected FD grids and frequencies, satisfactory results have been obtained with the use of the sparse LU decomposition of the system $A$ \cite{Operto2015,Kostin2019}.
In general, however, a parallel on-the-fly application of $A^{-1}$ and $A^{-*}$ is the only option.
The system matrix in \eqref{eq:KKT} is sparse complex-valued  Hermitian indefinite.
It is many times larger than the normal matrix, and it precludes from applying direct solvers in any form.
Its iterative solution is considered in the next section.

\section{Iterative solution of the KKT system}

\begin{figure}[h!]
	\centering
	\includegraphics[width=6cm]{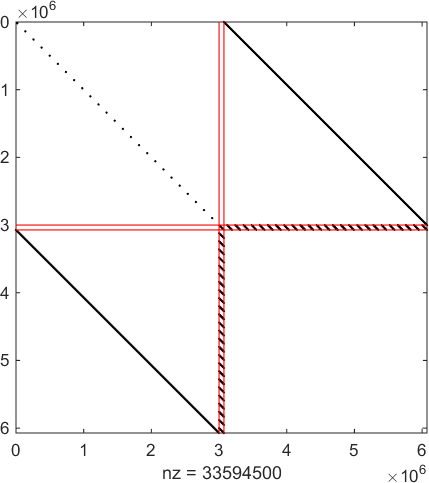}
	\caption{Partitioning (red lines) and non-zero elements (black dots) of the KKT matrix for a test problem. The problem contains 25 sources, each system matrix is 120,000 $\times$ 120,000. The size of the $L$ block is 72,000 $\times$ 72,000.}
	\label{fig:M}
\end{figure}

In this section, we design an iterative solver for the KKT system \eqref{eq:KKT}.
The system matrix $\mathcal{M}$ is can be explicitly partitioned across sources as follows,
\begin{equation}
	\label{eq:partitioned_system_matrix}
	\left[
	\begin{array}{@{}c|c|c@{}}
		\begin{matrix}
			F_0&        & \\
			   & \ddots & \\
			   &        & F_{K-1}\\
		\end{matrix} 
		&O
		&
		\begin{matrix}
			A^*_0&        & \\
			   & \ddots & \\
			   &        & A^*_{K-1}\\
		\end{matrix}
		\\		
		\hline
		O 
		&L
		&
		\begin{matrix}
			-P^*_0 \hdots -P^*_{K-1}
		\end{matrix}
		\\
		\hline
		\begin{matrix}
			A_0&        & \\
			   & \ddots & \\
			   &        & A_{K-1}\\
		\end{matrix}
		&
		\begin{matrix}
			-P_0   \\
			\vdots  \\
			-P_{K-1}
		\end{matrix}
		&
		O
	\end{array}
	\right]
	\left[
	\begin{array}{c}
		\mydelta{u}_0    \\
		\vdots        \\
		\mydelta{u}_{K-1}\\
		\hline
		\mydelta{s}        \\ 
		\hline
		\lambda_0	  \\
		\vdots        \\
		\lambda_{K-1}	  
	\end{array}
	\right]
	=
	\left[
	\begin{array}{c}
	Q_0^* W_0^T W_0 r_0 \\
	\vdots          \\
	Q_{K-1}^* W_{K-1}^T W_{K-1} r_{K-1} \\
	\hline
	-L s_n  \\
	\hline
	0\\
	\vdots \\
	0
	\end{array}
	\right]
	.
\end{equation}
Block $F$ is a positive semidefinite matrix, which is very sparse.
If receivers are located at the grid nodes, then $F$ is a diagonal matrix where the number of non-zero entries in each $F_j$ block equals the number of data points for the $j$-th source.
Block $A$ is a block-diagonal matrix of full rank whose spectrum lies in the left and right half-planes.
Block $L$ is positive-definite here but may become positive semidefinite if $L$ represents the Laplacian. 
In our case, $L=\varepsilon I$ and so its diagonal.
Block $P$ is a rectangular matrix consisting of sparse blocks. 
All submatrices of $P$ are diagonal; see definition \eqref{eq:P}.
The pattern of the matrix $\mathcal{M}$ for a test problem consisting of 25 sources 
is presented in Figure~\ref{fig:M}.

Our preconditioner is based on a block-triangular decomposition,
which is one of the most efficient strategies for preconditioning saddle-point problems.
In the context of the KKT matrix, this approach was considered in \cite{Haber2001,Biros2005}.  

The matrix $\mathcal{M}$ has singular (1,1) and (3,3) blocks, and so the 
direct block-LU decomposition is not possible.
Thus, we perform block-LU factorization of the permuted matrix as follows,
\begin{equation}
	\label{eq:permuted-block-LU}
		\left[
		\begin{array}{ccc}
			 A & O   & -P \\
			 F & A^* & O  \\
			 O &-P^* & L    \\
		\end{array}
		\right]
		=
		\left[
		\begin{array}{ccc}
			 I         & O   & O           \\
			 F A^{-1}  & I   & O           \\
			 O         & -P^* A^{-*}   & I \\
		\end{array}
		\right]
		\left[
		\begin{array}{ccc}
			 A  & O  & -P         \\
			 O  & A^*& F A^{-1} P \\
			 O  & O & H           \\
		\end{array}
		\right],
\end{equation}
where $H = L + P^* A^{-*} F A^{-1} P$.
The $U$-factor is an excellent preconditioner since applying its inverse to the initial (permuted) matrix results in a matrix whose spectrum consists of a single eigenvalue equal to 1. 
This approach would be equivalent to the reduced-space method and has no advantages because it requires a linear solve with the reduced Hessian $H$.

Let us approximate \eqref{eq:permuted-block-LU} by assuming $F \approx 0$,
which implies $H \approx L$.
Plugging $F=0, H=L$ into \eqref{eq:permuted-block-LU}, we get
\begin{equation}
	\label{eq:permuted-precond}
		\left[
		\begin{array}{ccc}
			 I         & O   & O           \\
			 O  & I   & O           \\
			 O         & -P^* A^{-*}   & I \\
		\end{array}
		\right]
		\left[
		\begin{array}{ccc}
			 A  & O  & -P\\
			 O  & A^*& O \\
			 O  & O & L  \\
		\end{array}
		\right]
		=
		\left[
		\begin{array}{ccc}
			 A  & O  & -P   \\
			 O  & A^*& O    \\
			 O  & -P^* & L  \\
		\end{array}
		\right].
\end{equation}
By reversing the permutation, we define the preconditioner to the original system matrix $\mathcal{M}$ as follows,
\begin{equation}
	\label{eq:precond}
		\mathcal{P} :=
		\left[
		\begin{array}{ccc}
			 O & O   &  A^* \\
			 O & L   & -P^* \\
			 A &-P   & O    \\
		\end{array}
		\right].
\end{equation}
The preconditioned matrix , $\mathcal{P}^{-1} \mathcal{M}$ is non-symmetric; therefore,
a non-symmetric Krylov solver is required to solve it.
We selected the GMRes and abbreviated the resulting method as FSGN-GMRes.

The preconditioner \eqref{eq:precond} has several appealing properties.
First, it requires three forward solves with matrices $A$, $A^*$, and $L$, e.i. almost the same work as needed for a single iteration of the CG with a reduced Hessian.
Second, the $A$ and $A^*$ both split into $K$ independent linear systems.
Thus, the solutions with blocks (1,3) and (3,1) can be performed in parallel.
Third, the preconditioner performs approximation $H \approx L$, but other approximations to $H$ are possible.
In particular, one can use the L-BFGS approximation in place of block $L$, see \cite{Biros2005}.

Now, we define another preconditioner, which is obtained from $\mathcal{P}$ by approximating $A$ by another matrix  $\tilde{A}$ as follows,
\begin{equation}
	\label{eq:precond_approx}
		\tilde{\mathcal{P}}:=
		\left[
		\begin{array}{ccc}
			 O & O   &  \tilde{A}^* \\
			 O & L   & -P^* \\
			 \tilde{A} &-P   & O    \\
		\end{array}
		\right].
\end{equation}
In this paper, we utilized the incomplete LU (ILU) factorization to approximate $A$.
This choice is motivated by the fact that the ILU is widely available and requires no tuning apart from selecting the fill-in level.

We will demonstrate below that the FSGN-GMRes with the $\tilde{\mathcal{P}}$ preconditioner is superior to the RSGN-CG.
Still, it is important to recognize that more efficient strategies to approximate $A$ are likely to exist.
In particular, the shifted-Laplacian preconditioner leveraged with the multigrid \cite{Erlangga2004} deserves a dedicated study.

%%%%%%%%%%%%%%%%%%%%%%%%%%%%%%%%%%%%%%%%%%%%%%%%%%%%%%%%%%%%%%%%%%%%%%%%%%%%%%%%%%%
%%%%%%%%%%%%%%%%%%%%%%%%%%%%%%%%%%%%%%%%%%%%%%%%%%%%%%%%%%%%%%%%%%%%%%%%%%%%%%%%%%%
%%%%%%%%%%%%%%%%%%%%%%%%%%%%%%%%%%%%%%%%%%%%%%%%%%%%%%%%%%%%%%%%%%%%%%%%%%%%%%%%%%%
%%%%%%%%%%%%%%%%%%%%%%%%%%%%%%%%%%%%%%%%%%%%%%%%%%%%%%%%%%%%%%%%%%%%%%%%%%%%%%%%%%%
\section{Numerical experiments}

This section is dedicated to numerical experiments. 
The optimization part of FWI was programmed in Python, whereas the forward and adjoint simulation problem is coded in C++.
The forward problem was solved directly with a sparse LU factorization.
All calculations were serial.

\begin{figure}[h!]
	\centering
	\includegraphics[width=10cm]{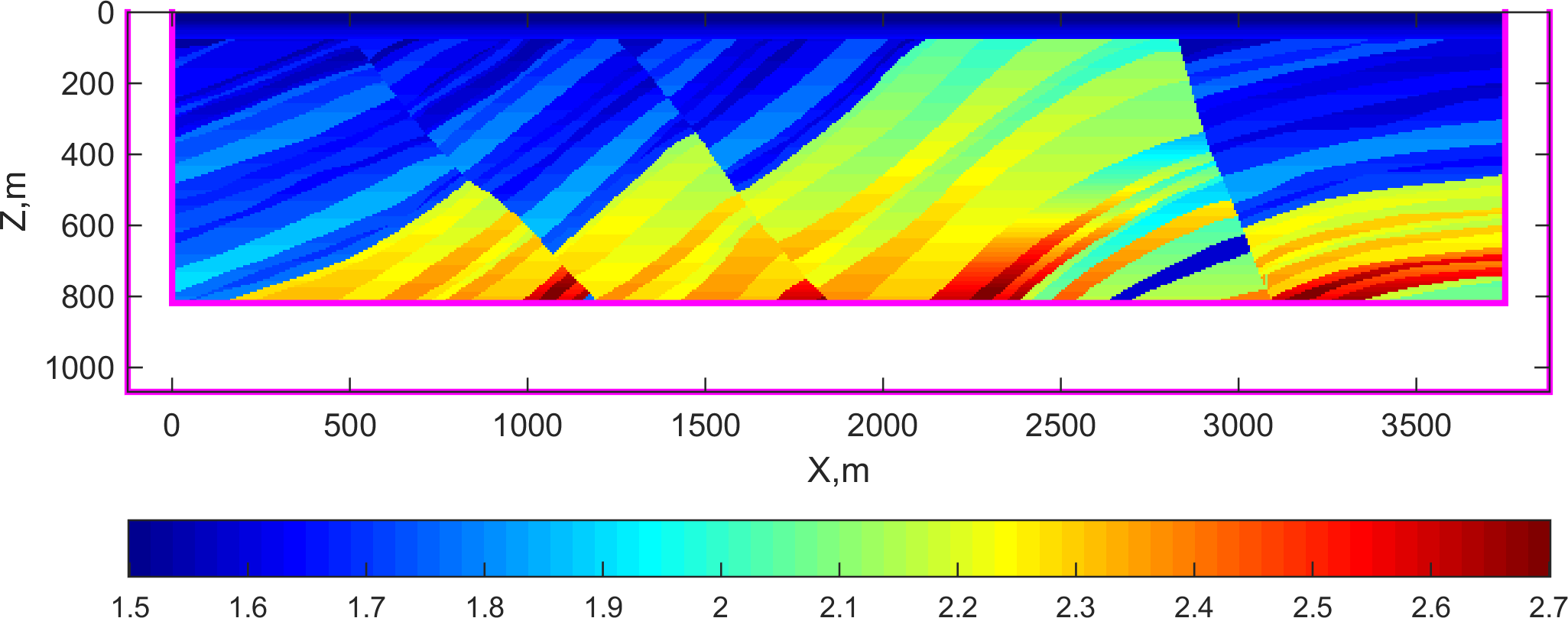}
	\caption{The distribution of the P-wave velocity, $c$,  used to simulate synthetic data. The numbers on the color scale are the velocity in km/s.
	The maroon lines on the sides indicate the boundaries of the PML region.}
	\label{fig:Model}
\end{figure}

We conducted numerical experiments using a part of the Marmousi velocity model (Figure~\ref{fig:Model}).
This P-wave velocity distribution was used in the data simulation, representing the geologic medium
we want to reconstruct.

\begin{figure}[H]
	\tiny
	\centering
	\begin{tabular}{c}
		\includegraphics[width=4in]{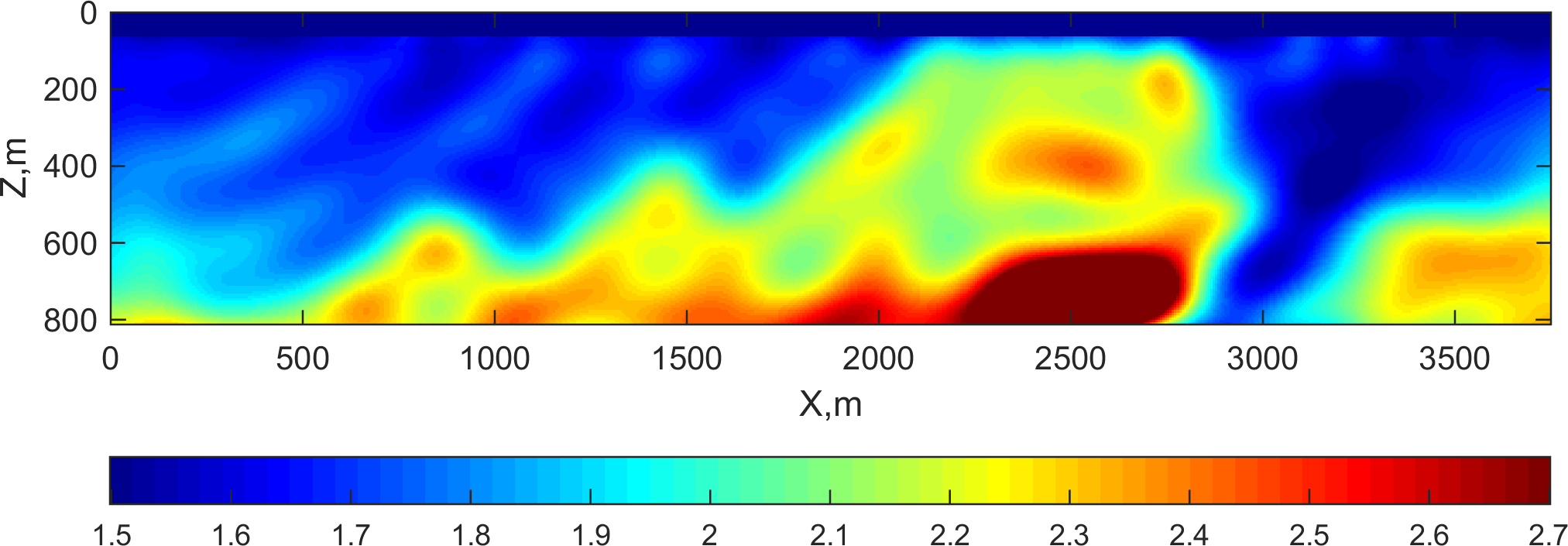}   \\ 
		(a) \\
		\includegraphics[width=4in]{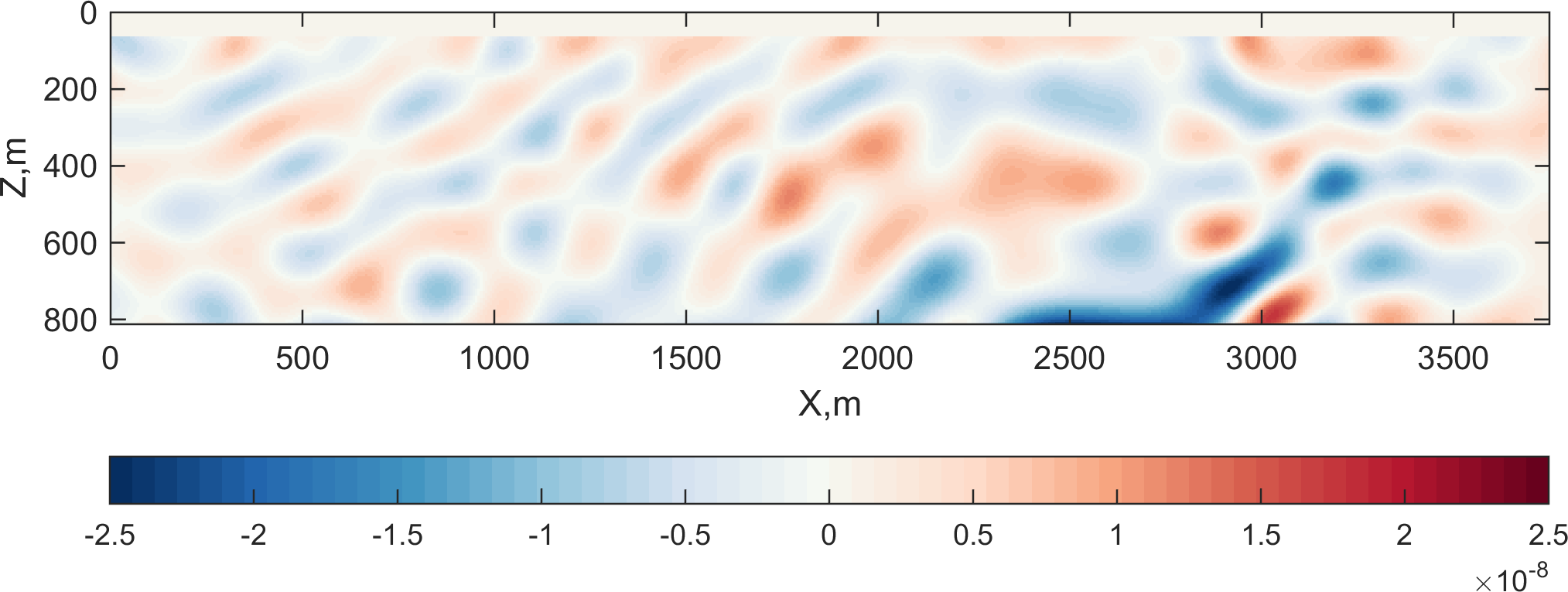}    \\
		(b) \\
		\includegraphics[width=4in]{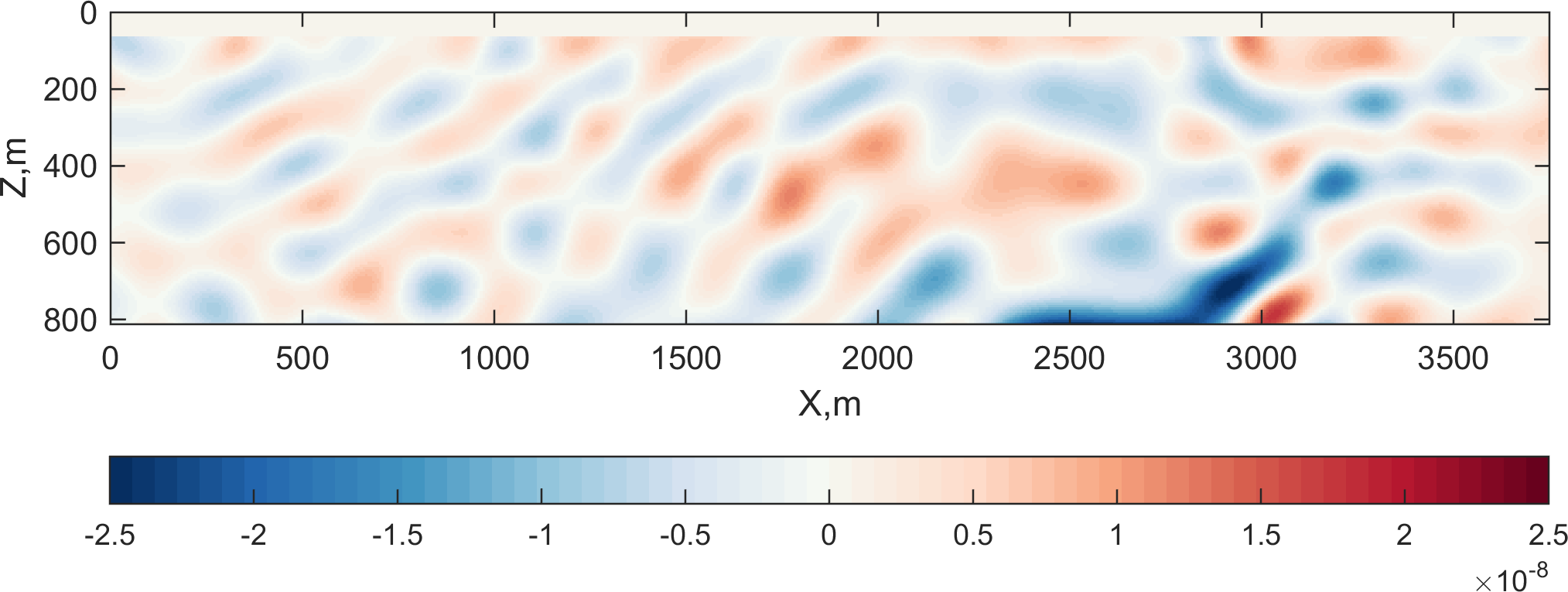} \\
		(c) \\
		\includegraphics[width=4in]{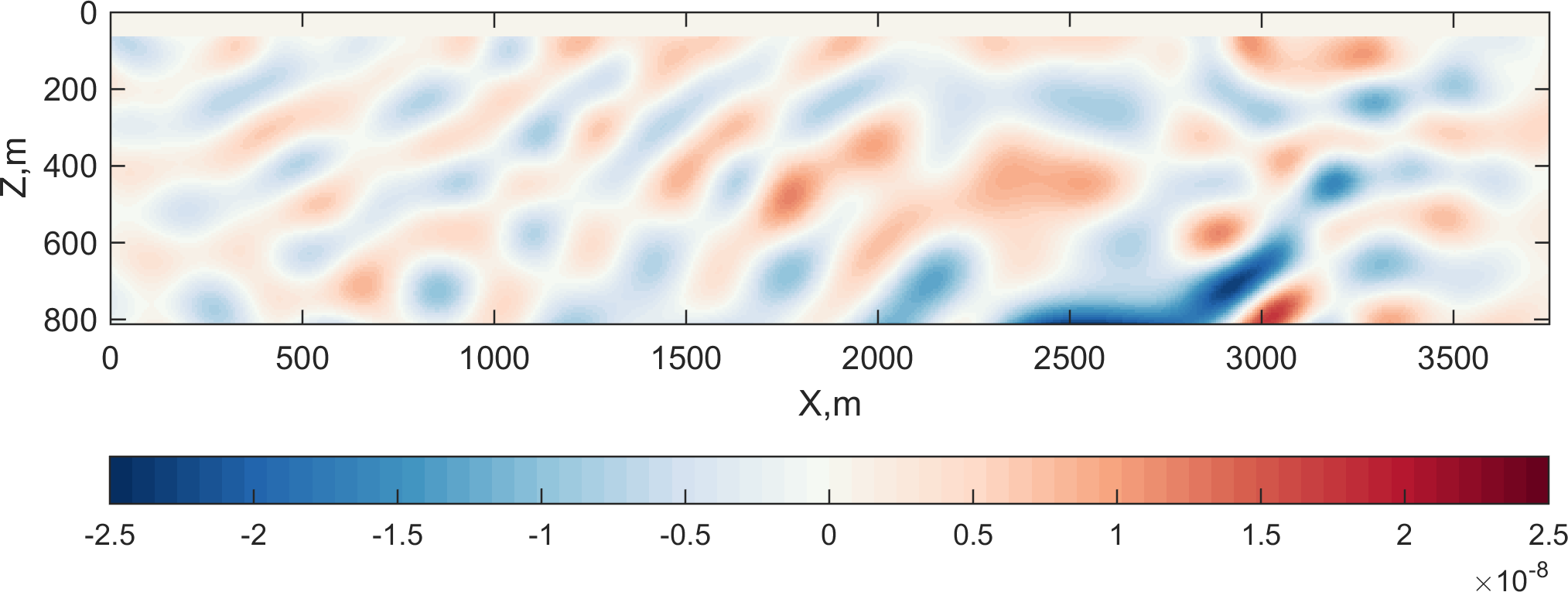} \\
		(d)\\		
	\end{tabular}	
		\caption{
		Comparison of the the two slowness updates computed at the end of a single GN iteration.
		Panel (a) shows the initial velocity distribution.
		Panel (b) depicts the slowness update obtained by applying the CG solver to the normal system of linear equations. 
		Panel (c) shows the slowness update by the GMRes solver with $\mathcal{P}$ applied to the KKT system.
		Panel (d) same as (c) but with preconditioner $\tilde{\mathcal{P}}$
		Note that in (a) color represents the P-wave velocity in km/s, 
		whereas in (b)-(d) color represents update to the squared slowness in $s^2/m^2$.
		}
		\label{fig:linear-solver-comparison}
\end{figure}

\begin{figure}[h!]
	\tiny
	\centering
	\includegraphics[width=3in]{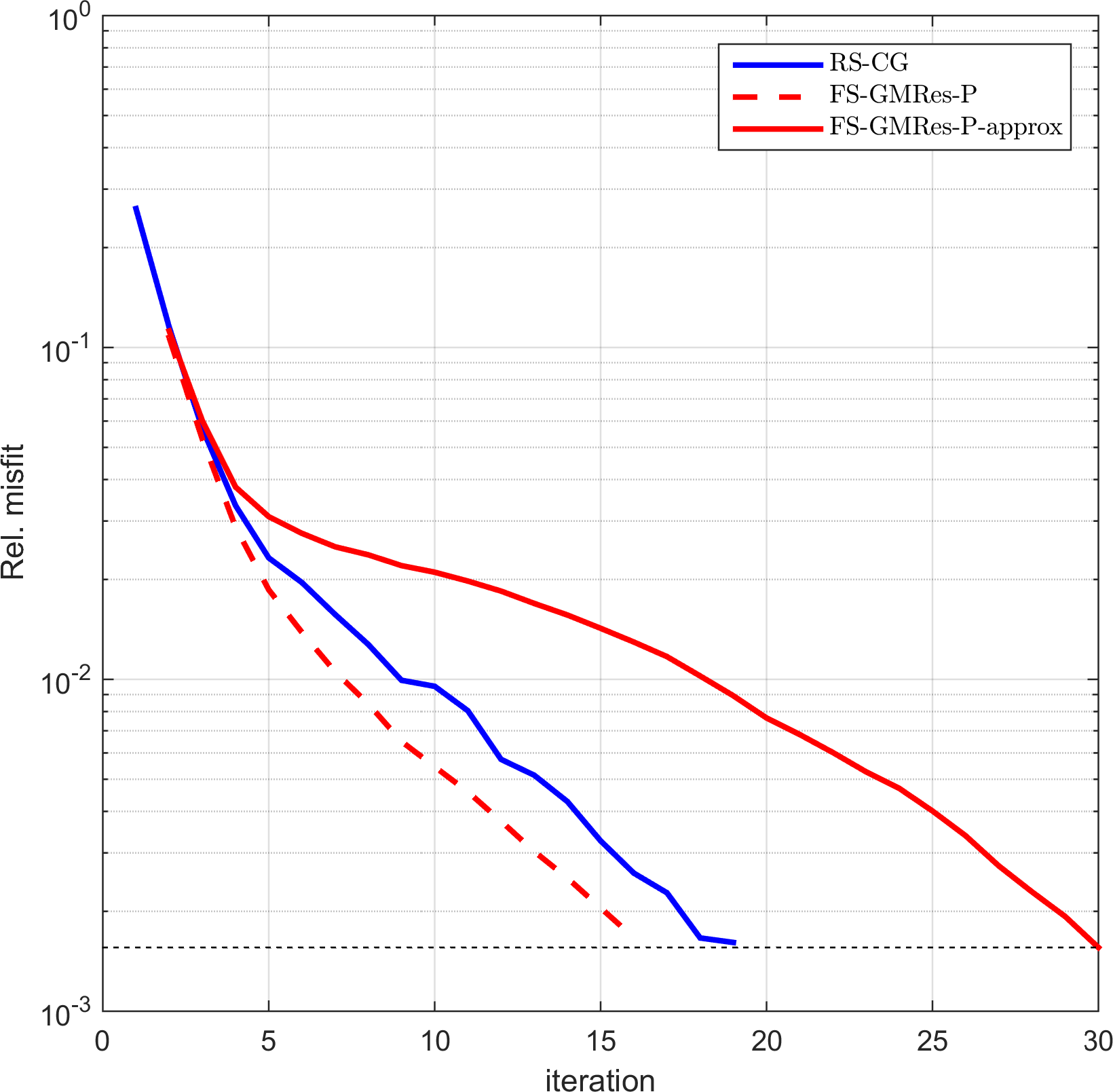}   \\ 
		\caption{Evolution of data misfit for the CG applied to the normal system and the preconditioned GMRes applied to the KKT matrix.
		The $Y$-axis represents the quantity $\norm{H \mydelta{s}_n-g}_2 / \norm{g}_2$.}
		\label{fig:linear-solver-convergence}
\end{figure}

The purpose of the first experiment was to compare, inside a single GN step, the performance of the CG applied to the normal system matrix with the preconditioned GMRes applied to the KKT matrix.
With the GMRes, we checked two preconditioners, $\mathcal{P}$ and $\tilde{\mathcal{P}}$.
The velocity distribution shown in Figure~\ref{fig:linear-solver-comparison}(a) was the initial model for the inversion.
The computed slowness updates are compared in Figure~\ref{fig:linear-solver-comparison}(b)-(d).
The preconditioner $\tilde{\mathcal{P}}$ was created by making use of ILU with the fill-in parameter 21.
It is not not possible to compare the CG and GMRes directly, because the CG computes $E_{cg} := \norm{H \mydelta{s}_n - g}_2$ at each iteration, 
whereas GMRes computes $E_{gmres} := \norm{\mathcal{M} \xi_n - b}_2$.
To make the comparison, we separately computed values of $E_{cg}$ from the GMres iterates $\xi_n$.
The convergence of the iterative methods are presented in Figure~\ref{fig:linear-solver-convergence}.
The CG and GMRes achieved approximately the same tolerance on iterations 19 and 30, respectively (Figure~\ref{fig:linear-solver-convergence}).
As expected, all methods generated almost identical updates.
The convergence rate of the CG is the same as that of the GMRes with the $\mathcal{P}$, 
although the GMRes happened to reach the desired tolerance by 3 iterations faster.
The convergence of the GMRes with $\tilde{\mathcal{P}}$ is slower than that with $\mathcal{P}$, but the cost is much lesser.
The performance of the algorithms is compared in Table~\ref{tab:1}.
We see the dramatic acceleration of matrix-vector multiplication in the GMRes compared to the CG.
This is because the linear solves with factorized matrices $\tilde{A}$ and $\tilde{A}^*$ inside $\tilde{\mathcal{P}}$ is much cheaper than the linear solves with $A$ and $A^*$ during the CG matrix-vector multiplication.
Another observation is that we needed a high value of fill-in (we used 21) to have an acceptable approximation $A \approx \tilde{A}$.
\begin{table}[h!]
    \centering
    \caption{Linear solvers comparison. The column 'CG' refers to the the CG applied to the normal system of linear equations, 
	'GMRes-full' is the GMRes applied to the KKT matrix preconditioned with $\mathcal{P}$, 
	'GMRes-approx' is the GMRes preconditioned with $\tilde{\mathcal{P}}$.}
    \begin{tabular}[t]{cccc}
    \hline
                           & CG    & GMRes-full &  GMRes-approx\\
    \hline
     Iterations            & 19    & 16         & 30     \\
     Time per iteration, s & 264   & 267        & 18     \\
	 ILU initialization, s & 0     & 0          & 121    \\
     Total time, s         & 4800  &  5391      & 834    \\
    \hline
    \end{tabular}
    \label{tab:1}
\end{table}%

\begin{figure}[h!]
	\tiny
	\centering
	\includegraphics[width=3in]{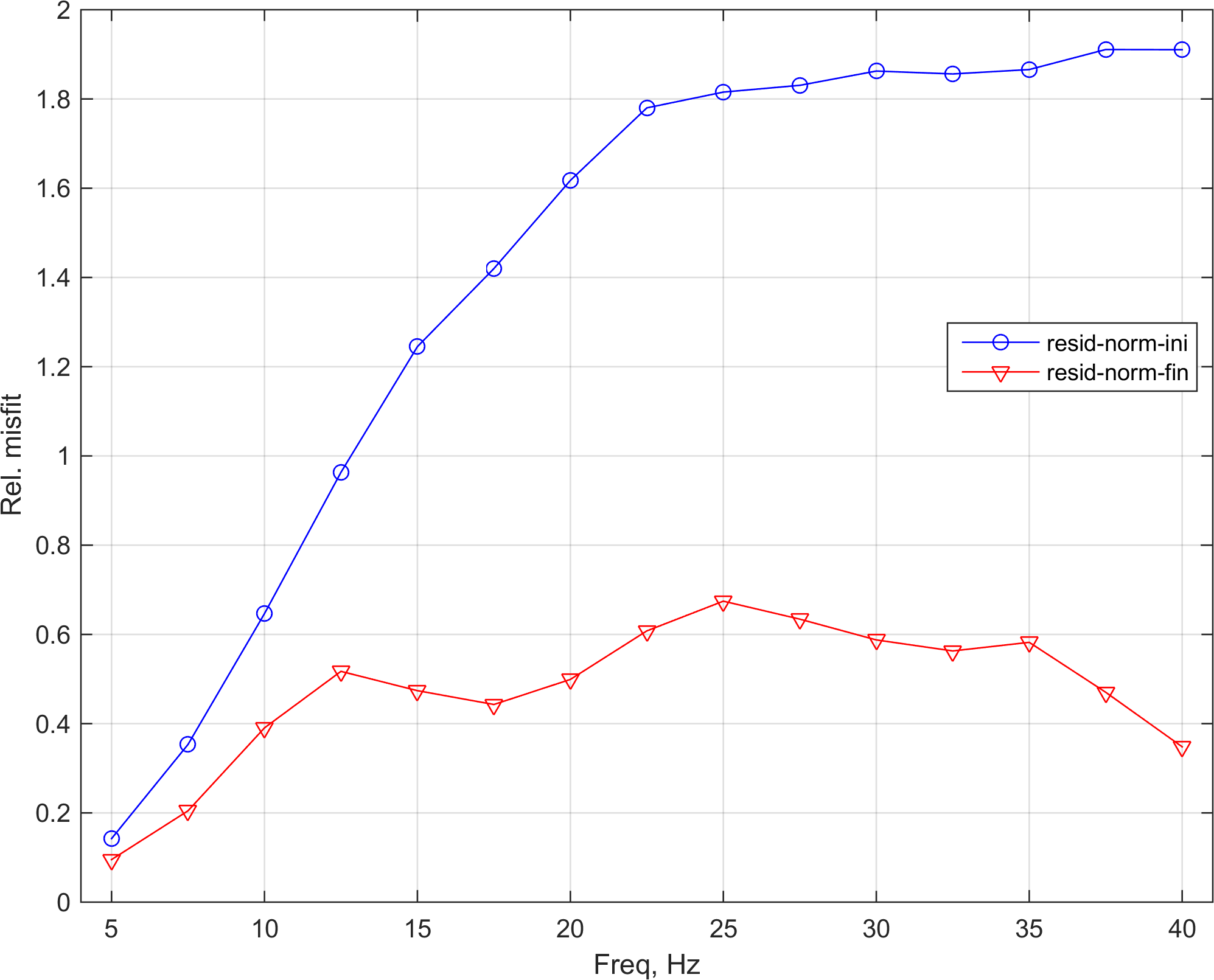}
		\caption{The relative residual misfit on Newton iterations:
		'resid-norm-ini' is $\norm{d_{initial}-d}_2/\norm{d}_2$,
		'resid-norm-fin' is $\norm{d_{final}-d}_2/\norm{d}_2$.
		}
		\label{fig:newton_misfit}
\end{figure}

In the second numerical experiment, we demonstrate that our approach serves as an efficient engine for the multi-frequency FWI.
The data were simulated at frequencies from 5 to 40~Hz with a step of 2.5~Hz.
At each frequency in the sequence, from lower to higher, 
a single GN step was performed, and the updated model was passed to the next frequency.
The regularization parameter $\varepsilon$ was changing linearly from 10  to 1E5.
The number of internal GMRes iterations was set to 30.
As indicated in Figure~\ref{fig:newton_misfit}, the data misfit was improved after the inversion.
The improvement becomes more pronounced at higher frequencies because the data magnitudes increase with the frequency and the initial model describes high-frequency data very poorly.

%With a high value of $\varepsilon$, the linear solver converges quickly at each Newton subproblem delivering a low tolerance. 
%However, the Newton subproblem itself becomes a poor approximation of the actual nonlinear problem, so the actual data fit deteriorates. 
%This trade-off between the amount of regularization and the quality of nonlinear iterations is well known.
%Thus, it is not characteristic of our approach.
%
\begin{figure}[H]
\small
\centering
%		\begin{tabular}{>{\centering} m{0.5in} >{\centering} m{0.9in} >{\centering} m{0.9in} m{0.9in} <{\centering} }
		\begin{tabular}{c}
			\includegraphics[width=4in]{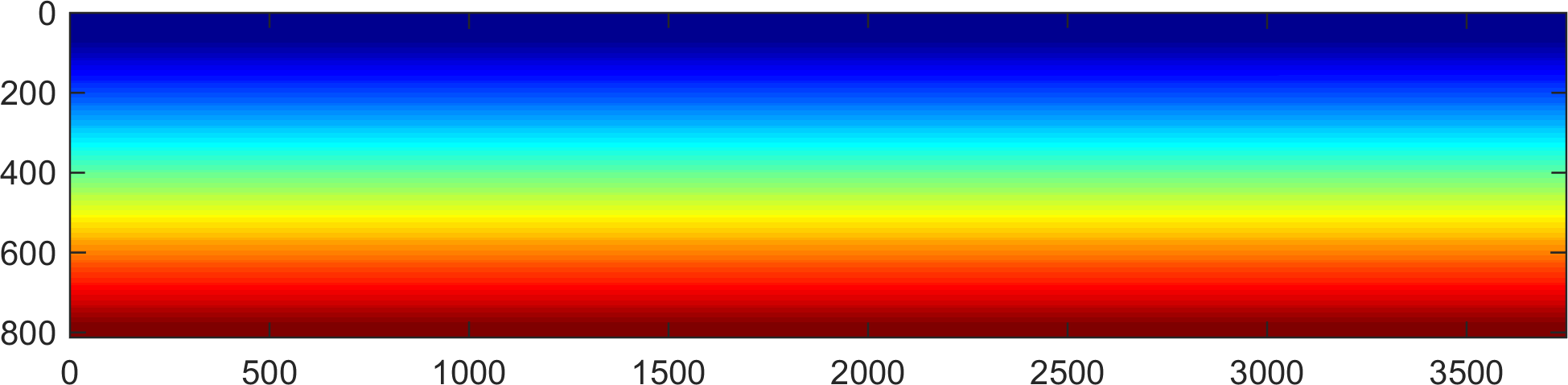} \\
			(a) Initial \\
			\includegraphics[width=4in]{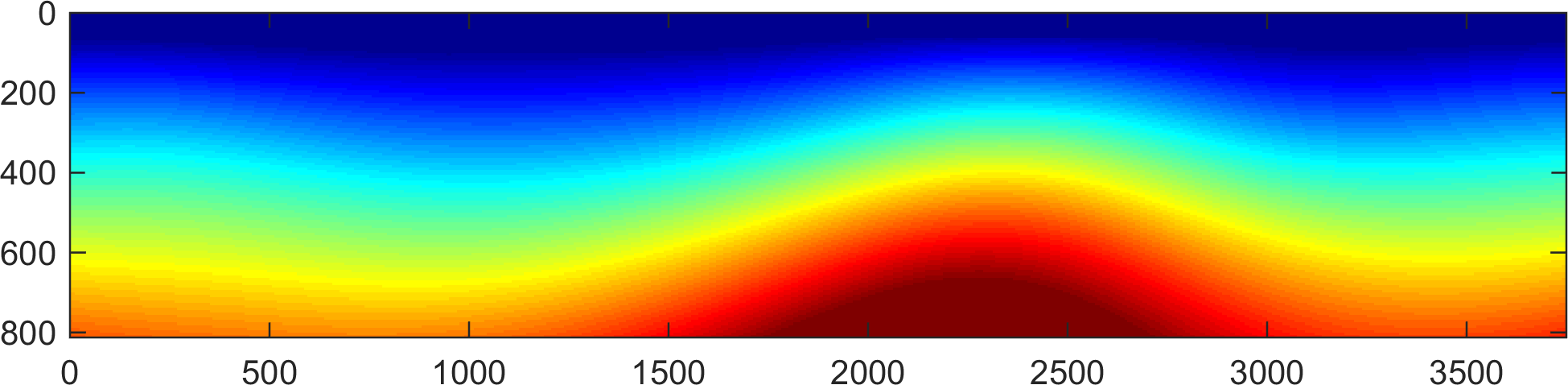}\\
			(b) 5~Hz \\
			\includegraphics[width=4in]{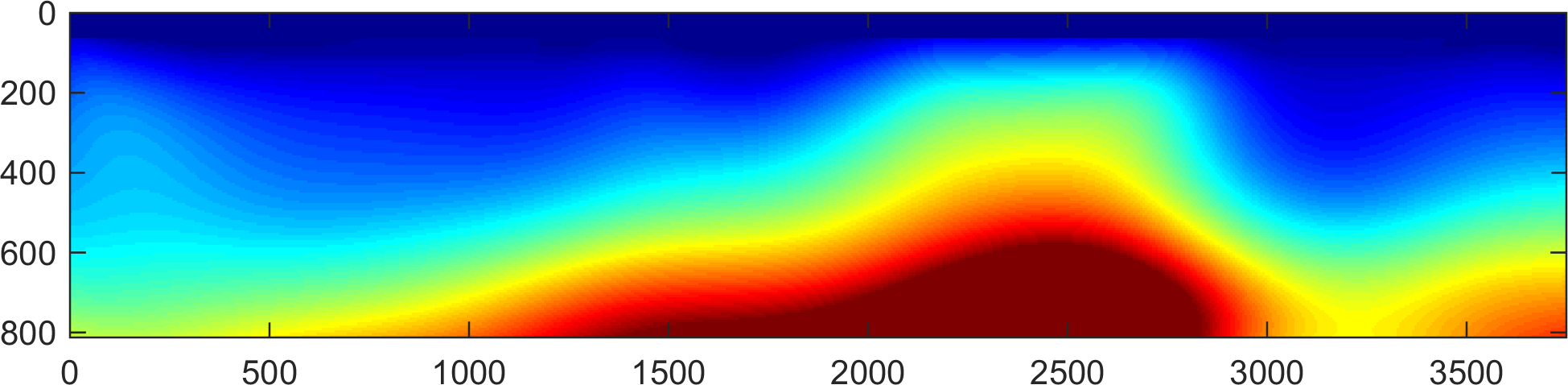} \\
			(c) 10~Hz \\
			\includegraphics[width=4in]{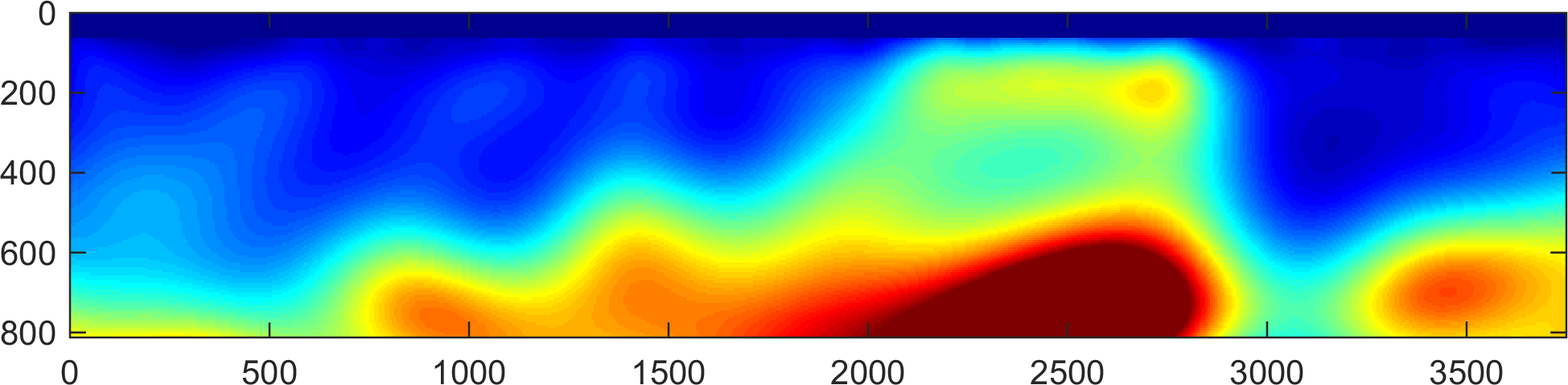} \\
			(d) 17.5~Hz \\
			\includegraphics[width=4in]{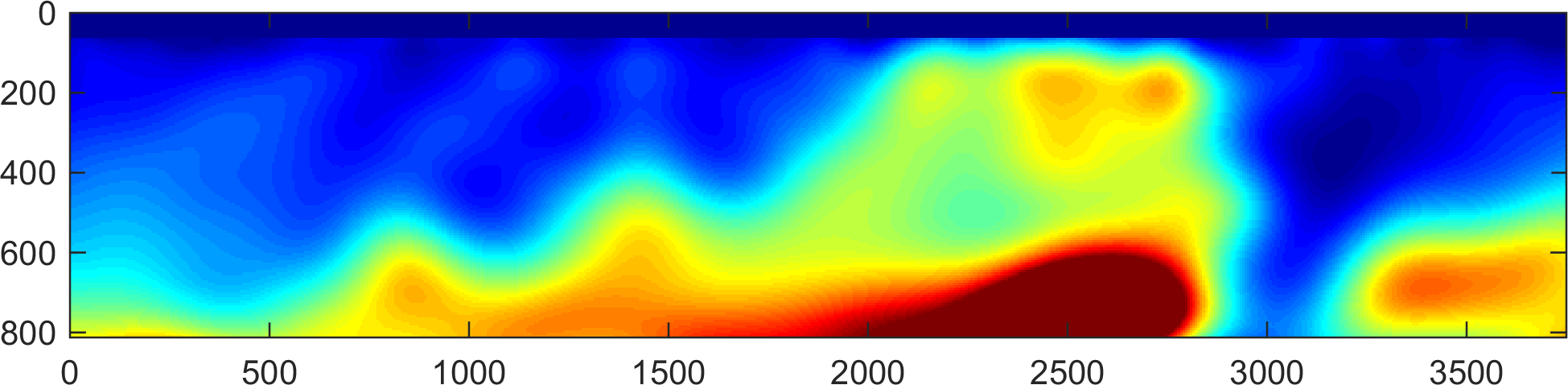}\\
			(e) 25~Hz \\
			\includegraphics[width=4in]{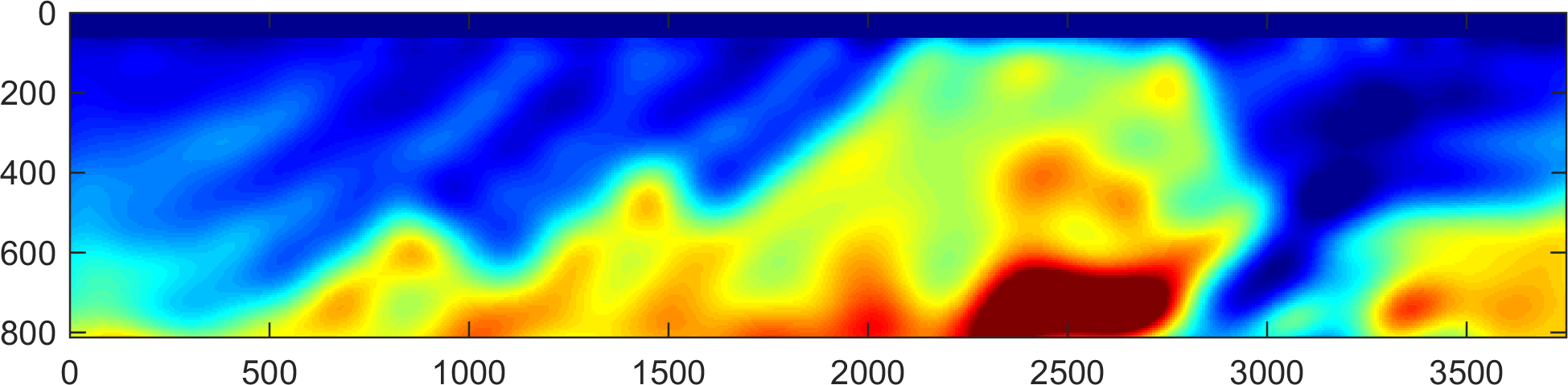}  \\
			(f) 40~Hz\\
		\end{tabular}
	\caption{The initial velocity distribution (a) and updated distributions (b)-(f) at different frequencies.	
	The color scales here and in Figure~\ref{fig:Model} are identical. The X- and Z-axis are marked in meters.
	Similar, if not identical, updates would have been obtained if we had performed the standard RSGN-CG optimization
	with the same accuracy of data fit.
	}
\label{fig:full-blown-fwi}
\end{figure}

The initial velocity model, several intermediate models, and the final model are presented in Figure~\ref{fig:full-blown-fwi}.
In this numerical experiment, we set the term $-L s_n$ in \eqref{eq:KKT} to zero and, thus, the algorithm was penalizing the norm of the update $\delta s$, not the updated model $s_n+\delta s$.
This was the reason for the high-velocity build-up at $X=2500$.
We stress again that identical FWI results would have been obtained if we had applied the standard RSGN-CG optimization (provided, of course, that the values of error $E_{cg}$ match in both approaches).
The difference between the two approaches is in their performance.

\section{Conclusions}

This paper proposes a new numerical method for the seismic FWI,
based on the full-space approach.
Prior work has demonstrated that the full-space approach can be used as a framework for efficient waveform inversion algorithms.
Those works have focused either on the penalty method with the primal-dual descend optimization
or the AL approach with the ADMM optimization procedure.
These approaches have demonstrated impressive results,
but may suffer from the slow convergence compared to the Newton iterations \cite{Boyd2010}.

This paper studies the iterative full-space FWI by the Newton method.
A preconditioned Krylov solver is applied to the KKT system on each Newton iteration. 
Our method computes the material parameter updates identical to those calculated by inverting the reduced Hessian, provided that the linear systems are solved exactly.
The critical element of our approach is a special preconditioner.
The preconditioned iterative algorithm requires two forward simulations (exact or approximate) at each iteration plus an inexpensive solve 
with the constraint matrix $L$.
We considered two variants of the preconditioner: 
preconditioner $\mathcal{P}$ computes linear solves with $A$ and $A^*$ exactly,
whereas preconditioner $\tilde{\mathcal{P}}$ computes the inexact solves.
A single application of $\mathcal{P}$ should 
therefore have almost the same cost compared to a single iteration of the CG applied to the normal system.
Our numerical experiments indicated that these two methods indeed have very close performance.
On the other hand, a single application of $\tilde{\mathcal{P}}$ is much cheaper.
In our numerical experiments, this approach reduced the execution time considerably.

For a real-scale 3D FWI, several problems need to be addressed. 
First, the forward and adjoint simulations need to be based on the theory of linear elasticity, so it remains to be shown that our approach works for elastic FWI.
Second, direct linear solvers become unfeasible in the 3D applications, even for the acoustic equations, 
not to mention the linear elasticity.
Thus, a reliable preconditioned iterative solver should be applied for the forward and adjoint simulations.
Substantial progress in this direction has been made; see \cite{Yavich2021} and references therein.
In principle, any preconditioner designed for the forward simulation can be used inside $\tilde{\mathcal{P}}$.
Third, the convergence can be improved by a more sophisticated approximation of the reduced Hessian,
such as the L-BFGS or the background Hessian.
Virtually any method designed to accelerate the reduced-space  Newton method can be used to compute the $H$-block solve within our approach.
Finally, using the BiCGStab instead of the GMRes will eliminate the memory needed to store GMRes iterates between restarts.
Future work, therefore, should include a study of these options.

\section{Acknowledgments}

This research was partially supported by the Russian Science Foundation, project no.
21-11-00139. The authors acknowledge computational resources granted by Complex
for Simulation and Data Processing for Mega-Science Facilities at NRC "Kurchatov
Institute", http://ckp.nrcki.ru/.

\bibliographystyle{unsrt}
\bibliography{fwi-oc}

\begin{thebibliography}{10}

\bibitem{Lailly1983}
P.~Lailly.
\newblock The seismic inverse problem as a sequence of before stack migrations.
\newblock In {\em Conference on Inverse Scattering, Theory and Application},
  pages 206--220. SIAM, 1983.

\bibitem{Tarantola1984}
A.~Tarantola.
\newblock Inversion of seismic reflection data in the acoustic approximation.
\newblock {\em Geophysics}, 49(8):1259--1266, 1983.

\bibitem{Mora1987}
P.~R. Mora.
\newblock Nonlinear two-dimensional elastic inversion of multi-offset seismic
  data.
\newblock {\em Geophysics}, 52(9):1211--1228, 1987.

\bibitem{Virieux2009}
J.~Virieux and S.~Operto.
\newblock An overview of full-waveform inversion in exploration geophysics.
\newblock {\em Geophysics}, 74:WCC1--WCC26, 2009.

\bibitem{Tromp2020}
J.~Tromp.
\newblock Seismic wavefield imaging of {E}arth’s interior across scales.
\newblock {\em Nature Reviews Earth \& Environment}, 1:40–53, 2020.

\bibitem{Plessix2009}
R.~Plessix.
\newblock Three-dimensional frequency-domain full-waveform inversion with an
  iterative solver.
\newblock {\em Geophysics}, 74(6):WCC149--WCC157, 2009.

\bibitem{Warner2013}
M.~Warner, A.~Ratcliffe, T.~Nangoo, J.~Morgan, A.~Umpleby, N.~Shah, V.~Vinje,
  I.~\v{S}tekl, L.~Guasch, C.~Win, G.~Conroy, and A.~Bertrand.
\newblock Anisotropic 3d full-waveform inversion.
\newblock {\em Geophysics}, 78(2):R59--R80, 2013.

\bibitem{Operto2015}
S.~Operto, A.~Miniussi, R.~Brossier, L.~Combe, L.~Métivier, V.~Monteiller,
  A.~Ribodetti, and J.~Virieux.
\newblock {Efficient 3-{D} frequency-domain mono-parameter full-waveform
  inversion of ocean-bottom cable data: application to {V}alhall in the
  visco-acoustic vertical transverse isotropic approximation}.
\newblock {\em Geophysical Journal International}, 202(2):1362--1391, 07 2015.

\bibitem{Operto2018}
S~Operto and A~Miniussi.
\newblock {On the role of density and attenuation in three-dimensional
  multiparameter viscoacoustic VTI frequency-domain FWI: an OBC case study from
  the North Sea}.
\newblock {\em Geophysical Journal International}, 213(3):2037--2059, 2018.

\bibitem{Fichtner2009}
A.~Fichtner, B.~L.~N. Kennett, H.~Igel, and H.-P. Bunge.
\newblock Full seismic waveform tomography for upper-mantle structure in the
  australasian region using adjoint methods.
\newblock {\em Geophysical Journal International}, 179(3):1703--1725, 2009.

\bibitem{Shin2001}
C.~Shin, S.~Jang, and D.-J. Min.
\newblock Improved amplitude preservation for prestack depth migration by
  inverse scattering theory.
\newblock {\em Geophysical Prospecting}, 49(5):592--606, 2001.

\bibitem{Mulder2008}
W.A. Mulder and R.-E. Plessix.
\newblock Exploring some issues in acoustic full waveform inversion.
\newblock {\em Geophysical Prospecting}, 56(6):827--841, 2008.

\bibitem{Brossier2009}
R.~Brossier, S.~Operto, and J.~Virieux.
\newblock Seismic imaging of complex onshore structures by 2d elastic
  frequency-domain full-waveform inversion.
\newblock {\em Geophysics}, 74(6):WCC105--WCC118, 2009.

\bibitem{Pratt1998}
R.~G. Pratt, C.~Shin, and G.~J. Hick.
\newblock {Gauss–Newton and full Newton methods in frequency–space seismic
  waveform inversion}.
\newblock {\em Geophysical Journal International}, 133(2):341--362, 1998.

\bibitem{Epanomeritakis2008}
I.~Epanomeritakis, V.~Ak\c{c}elik, O~Ghattas, and J~Bielak.
\newblock A {N}ewton-{CG} method for large-scale three-dimensional elastic
  full-waveform seismic inversion.
\newblock {\em Inverse Problems}, 24(3):034015, 2008.

\bibitem{Metivier2017}
L.~M{\'e}tivier, R.~Brossier, S.~Operto, , and J.~Virieux.
\newblock Full waveform inversion and the truncated newton method.
\newblock {\em SIAM Review}, 59(1):153–195, 2017.

\bibitem{Haber2000}
E.Haber, U.~M. Ascher, and D.~Oldenburg.
\newblock On optimization techniques for solving nonlinear inverse problems.
\newblock {\em Inverse Problems}, 16(4):1263, 2000.

\bibitem{Haber2001}
E.~Haber and U.~M. Ascher.
\newblock Preconditioned all-at-once methods for large, sparse parameter
  estimation problems.
\newblock {\em Inverse Problems}, 17:1847--1864, 2001.

\bibitem{Haber2004}
E.~Haber, U.~M. Ascher, and D.~W. Oldenburg.
\newblock Inversion of 3{D} electromagnetic data in frequency and time domain
  using an inexact all-at-once approach.
\newblock {\em Geophysics}, 69(5):1216--1228, 2004.

\bibitem{Abubakar2009}
A.~Abubakar, W.~Hu, T.~M. Habashy, and P.~M. van~den Berg.
\newblock Application of the finite-difference contrast-source inversion
  algorithm to seismic full-waveform data.
\newblock {\em Geophysics}, 74(6):WCC47--WCC58, 2009.

\bibitem{vanLeeuwen2013}
T.~van Leeuwen and F.~J. Herrmann.
\newblock {Mitigating local minima in full-waveform inversion by expanding the
  search space}.
\newblock {\em Geophysical Journal International}, 195(1):661--667, 2013.

\bibitem{vanLeeuwen2016}
T.~van Leeuwen and F.~J. Herrmann.
\newblock A penalty method for {PDE}-constrained optimization in inverse
  problems.
\newblock {\em Inverse Problems}, 32:015007, 2016.

\bibitem{Nocedal}
J.~Nocedal and S.~J. Wright.
\newblock {\em Numerical Optimization}.
\newblock Springer, New York, NY, USA, 2e edition, 2006.

\bibitem{Rizzuti2021}
G.~Rizzuti, M.~Louboutin, R.~Wang, and F.~J. Herrmann.
\newblock A dual formulation of wavefield reconstruction inversion for
  large-scale seismic inversion.
\newblock {\em Geophysics}, 86(6):R879--R893, 2021.

\bibitem{Boyd2010}
S.~Boyd, N.~Parikh, E.~Chu, B.~Peleato, and J.~Eckstein.
\newblock Distributed optimization and statistical learning via the alternating
  direction method of multipliers.
\newblock {\em Found. Trends Mach. Learn.}, 3(1):1–122, jan 2011.

\bibitem{Aghamiry2019}
H.~S. Aghamiry, A.~Gholami, and S.~Operto.
\newblock Improving full-waveform inversion by wavefield reconstruction with
  the alternating direction method of multipliers.
\newblock {\em Geophysics}, 84(1):R125--R148, 2019.

\bibitem{Aghazade2022}
K.~Aghazade, A.~Gholami, H.~S. Aghamiry, and S.~Operto.
\newblock Anderson-accelerated augmented lagrangian for extended waveform
  inversion.
\newblock {\em Geophysics}, 87(1):R79--R91, 2022.

\bibitem{Biros2005}
G.~Biros and O.~Ghattas.
\newblock Parallel {L}agrange--{N}ewton--{K}rylov--{S}chur methods for
  pde-constrained optimization. part i: The krylov--schur solver.
\newblock {\em SIAM Journal on Scientific Computing}, 27(2):687--713, 2005.

\bibitem{Yavich2021}
N.~Yavich, N.~Khokhlov, M.~Malovichko, and M.S Zhdanov.
\newblock Contraction operator transformation for the complex heterogeneous
  helmholtz equation.
\newblock {\em Computers and Mathematics with Applications}, 86:63--72, 2021.

\bibitem{Hustedt2004}
B.~Hustedt, S.~Operto, and J.~Virieux.
\newblock Mixed-grid and staggered-grid finite-difference methods for
  frequency-domain acoustic wave modelling.
\newblock {\em Geophysical Journal International}, 157(3):1269--1296, 2004.

\bibitem{Malovichko}
M.S. Malovichko, A.V. Tarasov, N.B. Yavich, and K.V. Titov.
\newblock Application of optimal control to inversion of self-potential data:
  theory and synthetic examples.
\newblock {\em IEEE Transactions on Geoscience and Remote Sensing}, in
  press:1--1, 2021.

\bibitem{Hinze2009}
M.~Hinze, R.~Pinnau, M.~Ulbrich, and S.~Ulbrich.
\newblock {\em Optimization with {PDE} Constraints}.
\newblock Springer, 2009.

\bibitem{Kostin2019}
V.~Kostin, S.~Solovyev, A.~Bakulin, and M.~Dmitriev.
\newblock Direct frequency-domain 3{D} acoustic solver with intermediate data
  compression benchmarked against time-domain modeling for full-waveform
  inversion applications.
\newblock {\em Geophysics}, 84(4):T207--T219, 2019.

\bibitem{Erlangga2004}
Y.A. Erlangga, C.~Vuik, and C.W. Oosterlee.
\newblock On a class of preconditioners for the helmholtz equation.
\newblock {\em Applied numerical mathematics}, (50):409--425, 2004.

\end{thebibliography}

\end{document}